\documentclass{amsart}

\DeclareMathSymbol{\twoheadrightarrow}  {\mathrel}{AMSa}{"10}

\def\Q{{\mathbf Q}}
\def\Z{{\mathbf Z}}

\def\F{{\mathbf F}}
\def\f{{\tilde F}}

\def\End{\mathrm{End}}
\def\Aut{\mathrm{Aut}}

\def\fchar{\mathrm{char}}

\def\tr{\mathrm{tr}}
\def\GL{\mathrm{GL}}
\def\SL{\mathrm{SL}}
\def\Sp{\mathrm{Sp}}
\def\M{\mathrm{M}}

\def\FF{F}
\def\f{f}
\def\g{h}

\def\OO{{\mathrm O}}
\def\O{{\mathcal O}}

\def\bmu{\boldsymbol \mu}
\newtheorem{thm}{Theorem}
\newtheorem{lem}[thm]{Lemma}

\newtheorem{prop}[thm]{Proposition}

\theoremstyle{definition}

\newtheorem{rem}[thm]{Remark}

\title[Representations of finite groups]
{Modular representations arising from self-dual $\ell$-adic
         representations of finite groups}

\author[A.\ Silverberg]{Alice Silverberg}
\address{Department of Mathematics, Ohio State University, 
Columbus, Ohio 43210, USA}
\email{silver\char`\@math.ohio-state.edu}

\author[Yu. G. Zarhin]{Yuri Zarhin}
\address{Department of Mathematics, Pennsylvania State University, 
University Park, PA 16802, USA}
\email{zarhin\char`\@math.psu.edu}

\thanks{Silverberg would like to thank NSA, NSF, and the Science Scholars
Fellowship Program at the Bunting Institute for financial support.
Zarhin would like to thank NSF for financial support.}

\begin{document}

\maketitle

Let $G$ be a finite subgroup of a symplectic group 
$\Sp_{2d}(\Q_\ell)$.
Despite the fact (\cite{symplectic}) 
that $G$ can fail to be conjugate
in $\GL_{2d}(\Q_\ell)$ to a subgroup of $\Sp_{2d}(\Z_\ell)$,
we prove that it can nevertheless be embedded in $\Sp_{2d}(\F_\ell)$
in such a way that the characteristic polynomials are preserved
(mod $\ell$), as long as $\ell>3$.

We start with the following ``rigidity" result,
which is in the spirit of similar results by Minkowski and Serre.

\begin{prop}
\label{minkthm}
Suppose $\ell$ is a prime number, and
$K$ is a discrete valuation field of characteristic zero
and residue characteristic $\ell$.
Let $\O$ denote the valuation ring and $\lambda$ the maximal
ideal. Let $e$ denote the ramification index of $K$ (i.e.,
$\ell\O=\lambda^e$), and suppose $2e<\ell-1$.
Suppose $S$ is a free $\O$-module of finite rank, $A$ is an
automorphism of $S$ of finite 
order, and $(A-1)^2\in\lambda\End(S)$.
Then $A=1$.
\end{prop}

\begin{proof}
This follows directly from Theorem 6.2 of \cite{Mink} with
$n=\ell$ and $k=2e$.
\end{proof}

Note that the hypothesis $2e<\ell-1$ is satisfied if
$e=1$ and $\ell \ge 5$. Next we state our main theorem. 

\begin{thm}
\label{mainthm}
Suppose 
$\ell$ is a prime number, and $K$ is a discrete valuation field 
of characteristic zero and residue characteristic $\ell$.
Let $\O$ denote the valuation ring, let $\lambda$ denote its maximal
ideal, and let $k$ denote the residue field $\O/\lambda$.
Let $e$ denote the ramification index of $K$, and 
suppose $2e<\ell-1$.
Suppose $V$ is a $K$-vector space of finite dimension $N$, 
suppose $\f:V \times V \to K$ 
is a nondegenerate alternating (respectively, symmetric) 
$K$-bilinear form, 
suppose $G$ is a finite group, and suppose 
$$\rho : G \hookrightarrow \Aut_K(V,{\f})$$
is a faithful representation of $G$ on $V$ that preserves
the form $f$.
Then there exist a nondegenerate alternating 
(respectively, symmetric) 
$k$-valued 
$k$-bilinear form 
${\f}_0$ on $k^N$, 
and a faithful representation
$${\bar \rho} : G \hookrightarrow \Aut_k(k^N,{\f}_0),$$
such that for every $g \in G$, the characteristic
polynomial of ${\bar \rho}(g)$ is the reduction modulo
$\lambda$ of the characteristic polynomial of $\rho(g)$.
\end{thm}

\begin{proof}
If $S$ is a $G$-stable $\O$-lattice in $V$, let
$$S^* = \{x \in V : {\f}(x,S) \subseteq \O \}.$$
Fix a $G$-stable $\O$-lattice $S$ in $V$.
Let $\pi$ denote a uniformizer for $\O$.
Multiplying ${\f}$ by an integral
power of $\pi$ if necessary, we may assume
that ${\f}(S,S) = \O$. Then $S \subseteq S^*$. 
Let $S_0=S$ and let
$$S_{i+1} = S_i + (\pi^{-1}S_i \cap \pi S_i^*) \quad 
\text{ for } \quad  i \ge 0.$$
Then $S_i$ is a $G$-stable $\O$-lattice in $V$,
${\f}(S_i,S_i) = \O$, and 
$S_i \subseteq S_{i+1} \subseteq S_{i+1}^* \subseteq S_i^*$.
Note that $S_{i+1} = S_i$ if and only if 
$\pi S_i^* \subseteq S_i$.
We have
$$S = S_0 \subseteq S_1 \subseteq S_2 \subseteq 
\ldots \subseteq S^*.$$
Since $S^*/S$ is finite, we have $S_j = S_{j+1}$
for some $j$. Let $T = S_j$.
Then $T$ is a $G$-stable 
$\O$-lattice in $V$ such that ${\f}(T,T) = \O$ and 
$\lambda T^* = \pi T^* \subseteq T$.

Let ${\bar \f} : T/\lambda T \times T/\lambda T \to k$
be the reduction of ${\f}$ modulo $\lambda$.
Then $\ker({\bar \f}) = \lambda T^*/\lambda T \cong T^*/T$.
Clearly, ${\bar \f}$ is nondegenerate on
$(T/\lambda T)/\ker({\bar {\f}}) \cong T/\lambda T^*$.
On $T^* \times T^*$, the form $\pi {\f}$ is $\O$-valued.
Let ${\tilde {\f}}$ denote the reduction modulo $\lambda$
of the restriction of $\pi {\f}$ to $T^* \times T^*$. 
Since $(T^*)^*=T$, we have $\ker({\tilde {\f}}) = T/\lambda T^*$. 
Therefore, ${\tilde {\f}}$ is nondegenerate
on $(T^*/\lambda T^*)/(T/\lambda T^*) \cong T^*/T \cong 
\ker({\bar {\f}})$.
We thus obtain a homomorphism
$$\psi : G \to 
\Aut_k((T/\lambda T)/\ker({\bar {\f}}),{\bar {\f}}) \times
\Aut_k(\ker({\bar {\f}}),{\tilde {\f}}) \cong$$
$$\Aut_k(T/\lambda T^*,{\bar {\f}}) \times
\Aut_k(T^*/T,{\tilde {\f}})
\hookrightarrow \Aut_k(T^*/\lambda T^*, {\bar {\f}} \times {\tilde {\f}}) 
\cong \Aut_k(k^N,{\f}_0)$$
for an appropriate pairing ${\f}_0$.
All the elements $\sigma\in\ker(\psi)$ act 
as the identity on $T/\lambda T^*$ and on $T^*/T$, and thus
$(\sigma-1)^2T^* \subseteq \lambda T^*$.
Proposition \ref{minkthm} implies $\psi$ is injective.  
Let ${\bar \rho} = \psi$.
\end{proof}

\begin{rem}
If one is concerned only with preserving the characteristic polynomials, 
and does not insist that ${\bar \rho}$ be an embedding, then in the
symplectic case the above
result can instead be achieved with the aid of Proposition 8 
of \cite{Serre1996}, rather than Proposition \ref{minkthm} above.
\end{rem}

\begin{thm}
\label{unitaryunram}
Suppose 
$\ell$ is a prime number, and $L$ is a discrete valuation field 
of characteristic zero and residue characteristic $\ell$.
Suppose $K$ is a quadratic extension of $L$, let $\O$ denote 
its valuation ring, let $\lambda$ denote the maximal
ideal, and let $k = \O/\lambda$. 
Let $e$ denote the ramification index of $K$, and 
suppose $2e<\ell-1$.
Suppose $V$ is a $K$-vector space of finite dimension $N$, 
suppose $\f:V \times V \to K$ is a nondegenerate 
pairing which is hermitian 
(respectively, skew-hermitian) with respect to the extension
$K/L$, suppose $G$ is a finite group, and 
suppose $\rho : G \hookrightarrow \Aut_K(V,{\f})$
is a faithful representation.
Suppose $K/L$ is unramified and thus 
$k$ is a quadratic extension of the residue
field $k_L$ of $L$. Then there exist a 
nondegenerate $k$-valued pairing ${\f}_0$ 
on $k^N$ which is hermitian  (respectively, skew-hermitian)
with respect to the
extension $k/k_L$, and 
a faithful representation
${\bar \rho} : G \hookrightarrow \Aut_k(k^N,{\f}_0)$,
such that for every $g \in G$, the characteristic
polynomial of ${\bar \rho}(g)$ is the reduction modulo
$\lambda$ of the characteristic polynomial of $\rho(g)$.
\end{thm}

\begin{proof}
Let $\pi$ denote a uniformizer for $L$. Then $\pi$ is
also a uniformizer for $K$, and 
the proof is a repetition of the proof of 
Theorem \ref{mainthm}. The reductions ${\bar \f}$ and 
${\tilde \f}$ are now hermitian (respectively, skew-hermitian).
\end{proof}

\begin{rem}
\label{unitaryram}
In the setting of Theorem \ref{unitaryunram},
suppose now that $K/L$ is ramified.
Then $k_L=k$, and one can write $K=L(\sqrt{D})$
where $D$ is a uniformizer for $L$ and
$\pi=\sqrt{D}$ is a uniformizer for $K$. Let $S$ be a
$G$-stable $\O$-lattice in $V$. Then 
$\pi^r\f(S,S)=\O$ for some integer $r$. If $r$ is even
then $\pi^r \f$ is hermitian (respectively,
skew-hermitian) and its reduction 
is symmetric (respectively, alternating). Further,
$\pi^{r+1}\f$ is skew-hermitian
(respectively, hermitian) and its reduction
is alternating (respectively, symmetric). 
If $r$ is odd, then 
$\pi^r\f$ is skew-hermitian (respectively,
hermitian) and its reduction is alternating 
(respectively, symmetric). Further, 
$\pi^{r+1}\f$ is hermitian
(respectively, skew-hermitian) and its reduction
is symmetric (respectively, alternating). 
In all cases one can proceed as in Theorem \ref{mainthm}
and obtain an embedding
of $G$ into a product of an orthogonal group $\OO_s(k)$ and a 
symplectic group $\Sp_{n-s}(k)$, which ``respects'' the 
characteristic polynomials. 
\end{rem}

\begin{lem}
\label{1101}
Suppose $\FF$ is a field of characteristic not equal to $2$ and
$f:\FF^2 \times \FF^2 \to \FF$
is a nondegenerate symmetric pairing. Let 
$g=\left(\begin{array}{cc}1&1\\0&1\end{array}\right)
\in \SL_2(\FF)$. Then $f$
is not $g$-invariant.
\end{lem}

\begin{proof}
Let $\{u,v\}$ denote the standard basis of $\FF^2$ over $\FF$.
Suppose $f$ is $g$-invariant. Then
$f(u,v)=f(gu,gv)=f(u,u+v)$, i.e.,
$f(u,u)=0$.
Also,
$f(v,v)=f(gv,gv)=f(u+v,u+v)=
f(v,v)+2f(u,v)$,
since $f$ is symmetric.
Since $\fchar(F) \ne 2$, we have 
$f(u,v)=0$.
Therefore $f(u,w)=0$ for every $w\in\FF^2$,
contradicting the nondegeneracy of $f$. 
\end{proof}

Let $Q_8$ denote the quaternion group of order $8$.
Let $\zeta_\ell$ denote a primitive $\ell$-th root of unity.
The next two results show that 
the condition $2e<\ell-1$ in Theorem \ref{mainthm} is sharp,
in both the symmetric and alternating cases.

\begin{prop}
\label{ramcounterex}
Let $\ell$ be an odd prime number, let
$K = \Q_\ell(\zeta_\ell+\zeta_\ell^{-1})$, 
let $V$ be $\Q_\ell(\zeta_\ell)$ viewed as a $2$-dimensional
$K$-vector space, and let $G=\bmu_\ell$.  
Then there exists a nondegenerate symmetric 
$K$-bilinear
form $\f:V \times V \to K$ such that there is a 
faithful irreducible representation  
$G \hookrightarrow \Aut_K(V,{\f})$. However, if
$\FF$ is a field of characteristic $\ell$ (in particular,
if $\FF=\F_\ell$), then
there does not exist a nondegenerate symmetric $\FF$-bilinear
form ${{\f}_0}:\FF^2\times\FF^2\to\FF$
such that $G$ embeds in $\Aut_\FF(\FF^2,{{\f}_0})$.
\end{prop}

\begin{proof}
Let $M=\Q_\ell(\zeta_\ell)$, and let $x \mapsto {\bar x}$ denote
the nontrivial automorphism of $M$ over $K$.
Define $f$ by
$f(x,y)=\tr_{M/K}(x{\bar y})$.
The desired injection is given by sending $\zeta_\ell$ to
multiplication by $\zeta_\ell$.
The nonexistence of ${\f}_0$ follows from Lemma \ref{1101},
since if $g$ is an element of $\GL_2(\FF)$ of order $\ell$, 
then there exists a basis of $\FF^2$ with respect 
to which $g=\left(\begin{array}{cc}1&1\\0&1\end{array}\right)$. 
\end{proof}

\begin{prop}
\label{fincounterex}
Let $\ell$ be an odd prime number, let
$K = \Q_\ell(\zeta_\ell+\zeta_\ell^{-1})$, 
let $V$ be $\Q_\ell(\zeta_\ell)^2$ viewed as a $4$-dimensional
$K$-vector space, and   
let $G = Q_8 \times \bmu_\ell$.
Then there exist a nondegenerate alternating 
$K$-bilinear form $\f:V \times V \to K$ 
and a faithful irreducible representation
$\rho : G \hookrightarrow \Aut_K(V,{\f})$ such that
there do not exist a field $\FF$ of characteristic $\ell$,
a nondegenerate alternating 
$\FF$-valued $\FF$-bilinear form ${{\f}_0}$ on $\FF^4$,
and a faithful representation
${\bar \rho} : G \hookrightarrow \Aut_\FF(\FF^4,\f_0)$,
having the property that for every $g \in G$, the characteristic
polynomial of ${\bar \rho}(g)$ is the reduction of 
the characteristic polynomial of $\rho(g)$.
\end{prop}

\begin{proof}
Let $M=\Q_\ell(\zeta_\ell)$, let $x \mapsto {\bar x}$ denote
the nontrivial automorphism of $M$ over $K$, and let $V_2$ be
$M$ viewed as a $2$-dimensional $K$-vector space.
Then $V=V_2 \times V_2 = \Q_\ell^2 \otimes_{\Q_\ell} V_2$.
Let $D_2$ denote the quaternion algebra over $\Q$ ramified
exactly at $2$ and infinity. Then $Q_8$ is a subgroup of
$D_2^\times$. The isomorphism 
$D_2 \otimes_\Q \Q_\ell \cong \M_2(\Q_\ell)$
induces a natural, absolutely irreducible, faithful representation  
$\rho_1 :  Q_8 \hookrightarrow \SL_2(\Q_\ell)
\cong \Aut(\Q_\ell^2,\f_1)$ where $\f_1$ is the standard
alternating pairing on $\Q_\ell^2$.
The natural action of $\bmu_\ell$ on $V_2$ 
respects the
nondegenerate symmetric bilinear form
${\f}_2 : V_2 \times V_2 \to K$ 
defined by ${\f}_2(x,y) = \tr_{M/K}(x{\bar y})$.
We thus obtain a faithful representation 
$\rho_2 : \bmu_\ell \hookrightarrow \Aut_K(V_2,\f_2)$.
Let $\rho$ be the representation
$$\rho : G = Q_8 \times \bmu_\ell \hookrightarrow 
\SL_2(\Q_\ell) \times \Aut_K(V_2,\f_2) \hookrightarrow \Aut_K(V,\f)$$
defined by $\rho(a,b) = \rho_1(a) \otimes \rho_2(b)$,
where $f=f_1\otimes f_2$. 
For every element $\sigma \in Q_8 \subset G$ of order $4$, 
the characteristic polynomial of $\rho(\sigma)$ is
$(t^2+1)^2$. 
For the element $-1 \in Q_8 \subset G$ of order $2$,
the characteristic polynomial of $\rho(-1)$ is
$(t+1)^4$. 

Suppose there exist a field $\FF$ of characteristic $\ell$,
a nondegenerate $\FF$-valued
alternating pairing ${\f_0}$ on $\FF^4$, and a
faithful representation
${\bar \rho} : G \hookrightarrow \Aut_\FF(\FF^4,\f_0)$
such that the characteristic polynomial  
of ${\bar \rho}(\sigma)$ is $(t^2+1)^2$ for every 
$\sigma\in Q_8$ of order $4$ and
the characteristic polynomial of ${\bar \rho}(-1)$ is
$(t+1)^4$. 
Let $V_0$ 
denote the corresponding faithful symplectic $F[G]$-module.
Since $\#Q_8$ is not divisible by $\ell$,
$V_0$ is a semisimple $F[Q_8]$-module.
By choosing a suitable basis we may assume that
$\rho_1(Q_8) \subset \SL_2(\Z_\ell)$.
The composition of $\rho_1$ with the reduction map gives 
a faithful representation
${\bar \rho}_1 : Q_8 \hookrightarrow \SL_2(\F_\ell)
\subseteq \SL_2(\FF)$.
The corresponding $\FF[Q_8]$-module $W$
is absolutely simple and symplectic.
By Schur's Lemma, every $Q_8$-invariant bilinear form
on $W$ is alternating. 
Since ${\bar \rho}_1 \oplus {\bar \rho}_1$ and the
restriction of ${\bar \rho}$ to $Q_8$ 
give rise to the same characteristic polynomials, 
the semisimple $F[Q_8]$-modules 
$V_0$ and $W \oplus W$ are isomorphic.
Since $\End_{Q_8}(W)=F$, we have $\End_{Q_8}(V_0)=M_2(F)$.
Fix a generator $c$ of $\bmu_\ell$. Then $c$ is an element
of $\End_{Q_8}(V_0)$ of multiplicative order $\ell$. 
We can therefore identify $V_0$ with 
$W \oplus W$ in such a way that
$c(x,y) = (x+y,y)$ for every $(x, y) \in W \oplus W = V_0$.
As in the proof of Lemma \ref{1101}, for $x, y \in W$ we have 
$${{\f}_0}((x,0),(0,y)) = {{\f}_0}(c(x,0),c(0,y)) = 
{{\f}_0}((x,0),(y,y)) = $$
$${{\f}_0}((x,0),(0,y)) + {{\f}_0}((x,0),(y,0)).$$
Therefore, 
\begin{equation}
\label{equat}
{{\f}_0}((x,0),(y,0)) = 0 \text{ for all }  x, y \in W.
\end{equation}
Further,
$${{\f}_0}((0,x),(0,y)) = {{\f}_0}(c(0,x),c(0,y)) = 
{{\f}_0}((x,x),(y,y)) = $$
$${{\f}_0}((0,x),(0,y)) + {{\f}_0}((x,0),(y,0)) 
+ {{\f}_0}((x,0),(0,y)) + {{\f}_0}((0,x),(y,0)).$$
Therefore,
${{\f}_0}((x,0),(0,y)) + {{\f}_0}((0,x),(y,0)) = 0$.
Since ${{\f}_0}$ is alternating, 
\begin{equation}
\label{equation1}
{{\f}_0}((x,0),(0,y)) = -{{\f}_0}((0,x),(y,0)) = {{\f}_0}((y,0),(0,x)).
\end{equation}

Define $\g : W \times W \to \FF$
by $\g(x,y)={{\f}_0}((x,0),(0,y))$.
By (\ref{equation1}), $\g$ is symmetric. 
Since ${{\f}_0}$ is $Q_8$-invariant, so is $\g$.
The nondegeneracy of $\g$ follows from 
(\ref{equat}) and the nondegeneracy of $f_0$.
Since $\g$ is a $Q_8$-invariant pairing on $W$, it is alternating.
Since $\g$ is both alternating and symmetric 
we have $\g=0$, giving a contradiction.
\end{proof}

\end{document}